\documentclass[onefignum,onetabnum]{siamart190516}

\usepackage{lipsum}
\usepackage{amsfonts}
\usepackage{graphicx}
\usepackage{epstopdf}
\usepackage{algorithmic}
\usepackage{graphicx}
\usepackage{color,latexsym}
\usepackage{verbatim}
\usepackage{lscape}
\usepackage{amsmath}
\usepackage{amssymb}
\usepackage{esint}
\usepackage{multirow}
\usepackage{hhline}
\usepackage{tabularx}
\usepackage{mathtools}
\usepackage{float}
\usepackage{setspace}

\newcommand{\defeq}{\vcentcolon=}

\ifpdf
  \DeclareGraphicsExtensions{.eps,.pdf,.png,.jpg}
\else
  \DeclareGraphicsExtensions{.eps}
\fi

\newsiamremark{remark}{Remark}
\newsiamremark{hypothesis}{Hypothesis}
\crefname{hypothesis}{Hypothesis}{Hypotheses}
\newsiamthm{claim}{Claim}

\headers{Approximation of Integrals Over the Volume of the Ball}{J. A. Reeger}

\title{Approximate Integrals Over the Volume of the Ball\thanks{Acknowledgement: This is a pre-print of an article published in the Journal of Scientific Computing. The final
authenticated version titled "Approximate Integrals Over the Volume of the Ball" is available online at: https://
doi.org/10.1007/s10915-020-01231-y.
\funding{This work was funded by the Office of Naval Research program Atmospheric Propagation Sciences for High Energy Lasers and the Air Force Office of Scientific Research project Radial Basis Functions for Numerical Simulation.}}}

\author{Jonah A. Reeger\thanks{
  (\email{jonah.reeger@gmail.com}, \url{}).}}

\usepackage{amsopn}

\begin{document}

\maketitle

\begin{abstract}
 A Radial Basis Function Generated Finite-Differences (RBF-FD) inspired technique for evaluating definite integrals over the volume of the ball in three dimensions is described.  Such methods are necessary in many areas of Applied Mathematics, Mathematical Physics and myriad other application areas.  Previous approaches needed restrictive uniformity in the node set, which the algorithm presented here does not require.  By using RBF-FD approach, the proposed algorithm computes quadrature weights for $N$ arbitrarily scattered nodes in only $O(N\mbox{ log}N)$ operations with high orders of accuracy.
\end{abstract}



\section{Introduction}
This article is concerned with the development of a method for the approximate evaluation of definite integrals over the volume of the ball in $\mathbb{R}^{3}$ of radius $\rho>0$.  That is, consider volume integration over the domain $\Omega=\left\{\mathbf{x}\in\mathbb{R}^{3}|\left\lVert\mathbf{x}-\mathbf{x}_{0}\right\rVert_{2}\leq\rho\right\}$, with $\mathbf{x}_{0}\in\mathbb{R}^{3}$.  Applications of this common problem in mathematics abound from the more specific scenarios of estimating the volumes of wells containing hydrocarbons beneath the Earth's surface \cite{PSMC2016} and recovering important diagnostic information in, e.g., Thermoacoustic or Photoacoustic Tomogrpahy \cite{PKLK2015} to some more generic problems in potential theory, magnetism and other concepts in mathematical physics \cite{JGL1922}.

The approximation of the values of definite integrals (quadrature when considering integration over an interval, or quadrature/cubature when considering integration over domains in two or more dimensions) is a rich topic dating back many centuries to early attempts to measure, for instance, the area of the circle \cite{WFMG2018}.  Since that time much research has been devoted to developing sophisticated and accurate techniques for estimating the values of integrals over intervals, areas, and volumes.  There are many texts devoted to summarizing such methods, see, for example \cite{AHS1971,ARKCWU1998,PKKMRS2005,WFMG2018}.

In the simplest case, the rules for quadrature over intervals are often constructed by replacing the integrand with a polynomial interpolant or a polynomial approximation and then integrating, or by enforcing that a rule be exact on a particular class of functions (like polynomials up to a particular degree).  This is how to arrive at, for instance,  Newton-Cotes or Gaussian Quadrature rules, respectively.  These rules for integration over intervals are then often leveraged in the context of iterated integrals by employing one-dimensional quadrature rules for each variable in turn, leading to so-called Cartesian product formulas.  Such formulas require structure in the node sets--spacing between nodes that is uniform or tied to the roots of orthogonal polynomials--across each variable of integration.  This requirement may be impractical or require an additional interpolation between the structured node set and locations where the integrand is specified.

To overcome these requirements on the structure of the node set, the value of the integral can be approximated utilizing concepts from number theory or through pseudorandom sampling of the integrand as in Monte-Carlo techniques \cite{ARKCWU1998}.  More common, however, is the alternative method of constructing quadrature rules over areas and volumes by replacing the integrand with a now multivariate polynomial interpolant or approximation and then integrating.  It is often the case that the basis set used for interpolation does not depend on the locations of the nodes (e.g. multivariate polynomials).  Such basis sets suffer from a question of the existence and uniqueness of an interpolant on unstructured node sets \cite{PCC1959,JCM1956}, which has prompted the use of Radial Basis Functions (RBFs) in the approximation of the integrand.

The proposed quadrature technique was developed to compute weights at whatever node locations are specified by the user.  This is because in applications, numerical quadrature is usually a follow-up to some other task (such as collecting data, or numerically solving PDEs), making it impractical to require node locations that are specific to the quadrature method. The present algorithm is therefore designed to find the quadrature weights given a node set defined by the application.  Further, the proposed algorithm allows for node sets featuring spatially varying separation when increased node density is needed to capture fine structure in the integrand.

The numerical method described in this paper is a generalization of RBF-FD (radial basis function-generated finite differences). This approach has so far mostly been used to approximate partial derivatives, with the key difference to regular finite differences that the node points no longer need to be grid-based (in particular, Cartesian node layouts are now known to be less-than-optimal \cite{NFGABLJW2016}. For surveys of RBFs and of RBF-FD methods as these are applied to PDEs, see \cite{BFNF2015a,BFNF2015b}.  Further, RBFs have been used successfully to construct quadrature rules for an interval in one-dimension, bounded domains in two-dimensions and, more specifically, integrals over bounded two-dimensional (piecewise-)smooth surfaces embedded in three-dimensions \cite{JARBF2016,JARBFMLW2016,JARBF2017}.

The following Section \ref{sec:Algorithm_Description} describes the present quadrature method.  Section \ref{sec:Test_Examples} describes some test examples, with illustrations of convergence rates and computational costs. Finally, section \ref{sec:Conclusions} outlines some conclusions.  A Matlab implementation of the method is available at Matlab Central's File Exchange \cite{BallVolumeQuadCodeMatlab}.

\section{Description of the key steps in the algorithm} \label{sec:Algorithm_Description}

Consider evaluating
\begin{align}
\iiint\limits_{\Omega}f(\mathbf{x})dV, \label{eq:volume_integral}
\end{align}
where $\Omega=\left\{\mathbf{x}\in\mathbb{R}^{3}:\left\lVert\mathbf{x}-\mathbf{x}_{0}\right\rVert_{2}\leq\rho\right\}$.

Similar to the work presented in \cite{JARBF2016,JARBFMLW2016,JARBF2017} for surface integrals, the proposed algorithm can be described in four steps:
\begin{enumerate}
    \item Decompose the domain of integration into $K\in \mathbb{Z}^{+}$ subdomains.
    \item On the $k^{\mbox{th}}$ subdomain construct an interpolant of the integrand.
    \item Integrate the interpolant of the integrand to determine weights for integrating a function $f$ over the $k^{\mbox{th}}$ subdomain.
    \item Combine the weights for the integrals over the $K$ subdomains to obtain a weight set for approximating the volume integral of $f$ over $\Omega$.
\end{enumerate}
Each of these steps is described in greater detail in what follows.

\subsection{Step 1: Decompose the domain of integration}

Consider first a set $\mathcal{S}_{N}=\left\{\mathbf{x}_{i}\right\}_{i=1}^{N}$ of $N$ points in $\Omega$, with a subset exactly on the boundary surface. On this set of points construct a tessellation $T=\left\{t_{k}\right\}_{k=1}^{K}$ (via Delaunay tessellation or some other algorithm) of $K$ tetrahedra.  These tetrahedra encompass the bulk of the volume of $\Omega$ but not the entire volume.

Let $\mathcal{K}_{S}\subset\{1,2,\ldots,K\}$ be the set of indices such that if $k\in \mathcal{K}_{S}$ then the tetrahedron $t_{k}$ has a face that is not shared with any of the other tetrahedra.  This face has all three vertices on the surface of $\Omega$ and, unless the surface is planar, there is a sliver of volume, $s_{k}$, between the unshared face and the spherical surface bounding the ball that must be accounted for in decomposing the volume.  Conversely, let $\mathcal{K}_{I}= \left\{1,2,\ldots,K\right\}\backslash \mathcal{K}_{S}$ be the indices of tetrahedra that do not have a face with three vertices on the surface.  With these definitions \eqref{eq:volume_integral} can be decomposed as
\begin{align}
\iiint\limits_{\Omega}f(\mathbf{x})dV=\sum\limits_{k\in \mathcal{K}_{I}}\iiint\limits_{t_{k}}f(\mathbf{x})dV+\sum\limits_{k\in \mathcal{K}_{S}}\left(\iiint\limits_{t_{k}}f(\mathbf{x})dV+\iiint\limits_{s_{k}}f(\mathbf{x})dV\right). \label{eq:volume_integral_decomposed}
\end{align}

\subsection{Step 2: Construct an interpolant of the integrand}

For each tetrahedron in $T$ define the sets $\mathcal{N}_{k}=\left\{\mathbf{x}_{k,j}\right\}_{j=1}^{n}$ to be the $n$ points in $\mathcal{S}_{N}$ nearest to the midpoint of $t_{k}$ (the average of the vertices of $t_{k}$).  Then for an individual $t_{k}$  the integral
\begin{align}
\iiint\limits_{t_{k}}f(\mathbf{x})dV
\end{align}
is evaluated by first approximating $f(\mathbf{x})$ by an RBF interpolant, with interpolation points from the set $\mathcal{N}_{k}$, and then integrating the interpolant. Often the RBF interpolant is a linear combination of (conditionally-) positive definite RBFs,
\begin{align}
\phi\left(\left\lVert \mathbf{x}-\mathbf{x}_{k,j}\right\rVert_{2}\right), j=1,2,\ldots,n\nonumber
\end{align}
and augmented by multivariate polynomial terms.  If $k\in\mathcal{K}_{S}$, then the same interpolant and set of nodes is used for approximating the integrand over $s_{k}$.  Define $\{\pi_{l}(\mathbf{x})\}_{l=1}^{M}$, with $M=\frac{(m+1)(m+2)(m+3)}{6}$, to be the set of all of the trivariate polynomial terms up to degree $m$. The interpolant is constructed as
\begin{align}
s(\mathbf{x}):=\sum_{j=1}^{n}c_{k,j}^{\mbox{RBF}}\phi\left(\left\lVert \mathbf{x}-\mathbf{x}_{k,j}\right\rVert_{2}\right)+\sum_{l=1}^{M}c_{k,l}^{p}\pi_{l}(\mathbf{x}),\nonumber
\end{align}
where $c_{k,1}^{RBF},\ldots,c_{k,n}^{RBF},c_{k,1}^{p},\ldots,c_{k,M}^{p}\in\mathbb{R}$ are chosen to satisfy the interpolation conditions $s(\mathbf{x}_{k,j})=f(\mathbf{x}_{k,j})$, $j=1,2,\ldots,n$, along with constraints $\sum_{j=1}^{n}c_{k,j}^{RBF}\pi_{l}(\mathbf{x}_{k,j})=0$, for $l=1,2,\ldots,M$.

\subsection{Step 3: Integrate the interpolant of the integrand}
By integrating the interpolant, the approximation of the integral of $f$ is reduced to
\begin{align}
\iiint\limits_{t_{k}}f(\mathbf{x})dV\approx\sum_{j=1}^{n}w_{k,j}f(\mathbf{x}_{k,j})\nonumber
\end{align}
for $k\in\mathcal{K}_{I}$ and
\begin{align}
\iiint\limits_{t_{k}}f(\mathbf{x})dV+\iiint\limits_{s_{k}}f(\mathbf{x})dV\approx\sum_{j=1}^{n}w_{k,j}f(\mathbf{x}_{k,j})\nonumber
\end{align}
for $k\in\mathcal{K}_{S}$.  A simple derivation can be carried out to show that the weights can be found by solving the linear system $A_{k}\mathbf{w}_{k}=\mathbf{I}_{k}$ with $(n+M)\times(n+M)$ matrix
\begin{align}
A_{k}=\left[\begin{array}{cc}\Phi_{k}^{T} & P_{k}\\P_{k}^{T} & 0\end{array}\right].\nonumber
\end{align}
The $n\times n$ submatrix $\Phi_{k}$ is made up of the RBFs evaluated at each point in $\mathcal{N}_{k}$, that is
\begin{align}
\Phi_{k,ij}=\phi\left(\left\lVert \mathbf{x}_{k,i}-\mathbf{x}_{k,j}\right\rVert\right),\mbox{ for }i,j=1,2,\ldots,n.\nonumber
\end{align}
Likewise the $n\times M$ matrix $P_{k,il}$ consists of the polynomial basis evaluated at each point in $\mathcal{N}_{k}$ so that
\begin{align}
P_{k,il}=\pi_{l}(\mathbf{x}_{k,i}),\mbox{ for }i=1,2,\ldots,n\mbox{ and }l=1,2,\ldots,M.\nonumber
\end{align}

If $k\in\mathcal{K}_{I}$, the right hand side, $\mathbf{I}_{k}$, includes integrals of the basis functions over $t_{k}$ only.  That is,
\begin{align}
I_{k,j}=\left\{\begin{array}{cc} \iiint\limits_{t_{k}}\phi\left(\left\lVert \mathbf{x}-\mathbf{x}_{k,j}\right\rVert\right)dV & j=1,2,\ldots,n \\ \iiint\limits_{t_{k}}\pi_{j-n}(\mathbf{x})dV & j=n+1,n+2,\ldots,n+M\end{array}\right..\nonumber
\end{align}
The integrals of the trivariate polynomial terms can be evaluated exactly via, for instance, the Divergence Theorem or barycentric coordinates.  For the RBFs, the integrals can be evaluated by further decomposing $t_{k}$ into four tetrahedra that share a common vertex.  Summing the integrals of the RBFs over the four tetrahedra results in the integral over $t_{k}$.  This process allows the volume integral over a tetrahedron to be reduced to four integrals in a single dimension.  Section \ref{sec:Tetrahedron_Integrals} explains this process.

On the other hand, for $k\in\mathcal{K}_{S}$
\begin{align}
I_{k,j}=\left\{\begin{array}{cc} \iiint\limits_{t_{k}}\phi\left(\left\lVert \mathbf{x}-\mathbf{x}_{k,j}\right\rVert\right)dV+\iiint\limits_{s_{k}}\phi\left(\left\lVert \mathbf{x}-\mathbf{x}_{k,j}\right\rVert\right)dV & j=1,2,\ldots,n \\ \iiint\limits_{t_{k}}\pi_{j-n}(\mathbf{x})dV+\iiint\limits_{s_{k}}\pi_{j-n}(\mathbf{x})dV  & j=n+1,\ldots,n+M\end{array}\right..\nonumber
\end{align}
The integrals over $t_{k}$ are evaluated using the methods described in the previous paragraph while the integrals over $s_{k}$ are approximated using a scheme discussed in section \ref{sec:Sliver_Integrals}.

\subsubsection{Integrals of RBFs Over Tetrahedra} \label{sec:Tetrahedron_Integrals}

Suppose that the tetrahedron $t_{k}$ has vertices $\mathbf{a}_{k}$, $\mathbf{b}_{k}$, $\mathbf{c}_{k}$ and $\mathbf{d}_{k}$, all points in $\mathbb{R}^{3}$.  Let $\mathbf{x}_{k,j}$ be some point in $\mathbb{R}^{3}$, which will be common to the four tetrahedra that will be integrated over to obtain the integral over $t_{k}$.  Although what follows applies for any point in $\mathbb{R}^{3}$, the point $\mathbf{x}_{k,j}$ is an interpolation node from the set $\mathcal{N}_{k}$ in this context.   A unit length normal vector to the side of $t_{k}$ with vertices $\mathbf{a}_{k}$, $\mathbf{b}_{k}$ and $\mathbf{c}_{k}$ is defined by
\begin{align}
\mathbf{n}_{\mathbf{a}_{k}\mathbf{b}_{k}\mathbf{c}_{k}}\defeq \frac{(\mathbf{b}_{k}-\mathbf{a}_{k})\times(\mathbf{c}_{k}-\mathbf{a}_{k})}{\left\lVert (\mathbf{b}_{k}-\mathbf{a}_{k})\times(\mathbf{c}_{k}-\mathbf{a}_{k})\right\rVert_{2}},\nonumber
\end{align}
and when defining a normal vector in what follows, the order of the vertices matters and should be taken as the order shown in this definition.  Further, let $\mathbf{e}_{k,j}$, $\mathbf{f}_{k,j}$, $\mathbf{g}_{k,j}$ and $\mathbf{h}_{k,j}$ be the orthogonal projections of $\mathbf{x}_{k,j}$ onto the sides of $t_{k}$ with vertices $\mathbf{a}_{k}$, $\mathbf{b}_{k}$ and $\mathbf{c}_{k}$; $\mathbf{a}_{k}$, $\mathbf{d}_{k}$ and $\mathbf{b}_{k}$; $\mathbf{a}_{k}$, $\mathbf{c}_{k}$ and $\mathbf{d}_{k}$; and $\mathbf{b}_{k}$, $\mathbf{d}_{k}$ and $\mathbf{c}_{k}$, respectively. For instance,
\begin{align}
\mathbf{e}_{k,j}=\mathbf{x}_{k,j}+\left[(\mathbf{a}_{k}-\mathbf{x}_{k,j})\cdot \mathbf{n}_{\mathbf{a}_{k}\mathbf{b}_{k}\mathbf{c}_{k}}\right]\mathbf{n}_{\mathbf{a}_{k}\mathbf{b}_{k}\mathbf{c}_{k}},
\end{align}
Then by applying the divergence theorem it can be shown that
\begin{align}
\iiint\limits_{t_{k}}\phi\left(\left\lVert \mathbf{x}-\mathbf{x}_{k,j}\right\rVert_{2}\right)dV&=
\nonumber\\
&\Bigg\{\mbox{sign}\left(\left(\mathbf{x}_{k,j}-\mathbf{e}_{k,j}\right)\cdot\mathbf{n}_{\mathbf{a}_{k}\mathbf{b}_{k}\mathbf{c}_{k}}\right)\iiint\limits_{t_{\mathbf{x}_{k,j}\mathbf{a}_{k}\mathbf{b}_{k}\mathbf{c}_{k}}}\phi\left(\left\lVert \mathbf{x}-\mathbf{x}_{k,j}\right\rVert_{2}\right)+\cdots\nonumber\\
&\mbox{sign}\left(\left(\mathbf{x}_{k,j}-\mathbf{f}_{k,j}\right)\cdot\mathbf{n}_{\mathbf{a}_{k}\mathbf{d}_{k}\mathbf{b}_{k}}\right)\iiint\limits_{t_{\mathbf{x}_{k,j}\mathbf{a}_{k}\mathbf{d}_{k}\mathbf{b}_{k}}}\phi\left(\left\lVert \mathbf{x}-\mathbf{x}_{k,j}\right\rVert_{2}\right)+\cdots\nonumber \\
&\mbox{sign}\left(\left(\mathbf{x}_{k,j}-\mathbf{g}_{k,j}\right)\cdot\mathbf{n}_{\mathbf{a}_{k}\mathbf{c}_{k}\mathbf{d}_{k}}\right)\iiint\limits_{t_{\mathbf{x}_{k,j}\mathbf{a}_{k}\mathbf{c}_{k}\mathbf{d}_{k}}}\phi\left(\left\lVert \mathbf{x}-\mathbf{x}_{k,j}\right\rVert_{2}\right)+\cdots\nonumber\\
&\mbox{sign}\left(\left(\mathbf{x}_{k,j}-\mathbf{h}_{k,j}\right)\cdot\mathbf{n}_{\mathbf{b}_{k}\mathbf{d}_{k}\mathbf{c}_{k}}\right)\iiint\limits_{t_{\mathbf{x}_{k,j}\mathbf{b}_{k}\mathbf{d}_{k}\mathbf{c}_{k}}}\phi\left(\left\lVert \mathbf{x}-\mathbf{x}_{k,j}\right\rVert_{2}\right)\Bigg\}.\nonumber
\end{align}

This expression for the integral over $t_{k}$ contains integrals over the four tetrahedra that share $\mathbf{x}_{k,j}$ as a common vertex.  Consider
\begin{align}
    \iiint\limits_{t_{\mathbf{x}_{k,j}\mathbf{a}_{k}\mathbf{b}_{k}\mathbf{c}_{k}}}\phi\left(\left\lVert \mathbf{x}-\mathbf{x}_{k,j}\right\rVert_{2}\right)dV\nonumber
\end{align}
since the remaining integrals over the tetrahedra are analogous.  The integrand is radially symmetric about $\mathbf{x}_{k,j}$ and depends only on the distance from $\mathbf{x}_{k,j}$ suggesting the change of variables
\begin{align}
    \mathbf{x}(\sigma,\lambda_{1},\lambda_{2})=\mathbf{x}_{k,j}+\sigma(\lambda_{1}\mathbf{a}_{k}+\lambda_{2}\mathbf{b}_{k}+(1-\lambda_{1}-\lambda_{2})\mathbf{c}_{k}-\mathbf{x}_{k,j}).\nonumber
\end{align}
In this change of variables, triangles similar to the side of $t_{k}$ with vertices $\mathbf{a}_{k}$, $\mathbf{b}_{k}$ and $\mathbf{c}_{k}$ are parameterized using barycentric coordinates and scaled by the nonnegative parameter $\sigma$, which also accounts for the distance of the triangle from the point $\mathbf{x}_{k,j}$.  Under this change of variables the integral becomes
\begin{align}
    \int\limits_{0}^{1}&\int\limits_{0}^{1-\lambda_{1}}\int\limits_{0}^{\sigma_{\mathbf{a}_{k}\mathbf{b}_{k}\mathbf{c}_{k}}}\phi\left(\sigma\left\lVert \lambda_{1}\mathbf{a}_{k}+\lambda_{2}\mathbf{b}_{k}+(1-\lambda_{1}-\lambda_{2})\mathbf{c}_{k}-\mathbf{x}_{k,j}\right\rVert_{2}\right)\sigma^{2}V_{k,j}{d\sigma}{d\lambda_{2}}{d\lambda_{1}},\nonumber
\end{align}
where $V_{k,j}=\left|(\mathbf{a}_{k}-\mathbf{x}_{k,j})\cdot\left[(\mathbf{b}_{k}-\mathbf{x}_{k,j})\times(\mathbf{c}_{k}-\mathbf{x}_{k,j})\right]\right|$ is six times the volume of the tetrahedron $t_{\mathbf{x}_{k,j}\mathbf{a}_{k}\mathbf{b}_{k}\mathbf{c}_{k}}$.  This volume appears in the Jacobian determinant from the change of variables.  Also,
\begin{align}
\sigma_{\mathbf{a}_{k}\mathbf{b}_{k}\mathbf{c}_{k}}=\frac{\left(\frac{1}{3}(\mathbf{a}_{k}+\mathbf{b}_{k}+\mathbf{c}_{k})-\mathbf{x}_{k,j}\right)\cdot\mathbf{n}_{\mathbf{a}_{k}\mathbf{b}_{k}\mathbf{c}_{k}}}{\left(\mathbf{c}_{k}-\mathbf{x}_{k,j}\right)\cdot \mathbf{n}_{\mathbf{a}_{k}\mathbf{b}_{k}\mathbf{c}_{k}}}\nonumber
\end{align}
is the value of $\sigma$ corresponding to the side of $t_{k}$ with vertices $\mathbf{a}_{k}$, $\mathbf{b}_{k}$ and $\mathbf{c}_{k}$.


Now, in the case of $\phi(r)=r^{2p+1}$, $p=0,1,2,\ldots$, the iterated integrals in $\sigma$ and then $\lambda_{2}$ can be computed in closed form. 
However, exploring the integration over $\lambda_{1}$ in Mathematica indicates the cost of a closed form expression for the integral is computationally too expensive, so the proposed algorithm uses standard pseudospectral methods for evaluating the integrals over $\lambda_{1}$.

\subsubsection{Integrals Over Slivers of Volume at the Surface} \label{sec:Sliver_Integrals}

\begin{figure}[h]
\begin{center}
\includegraphics[width=2.5in]{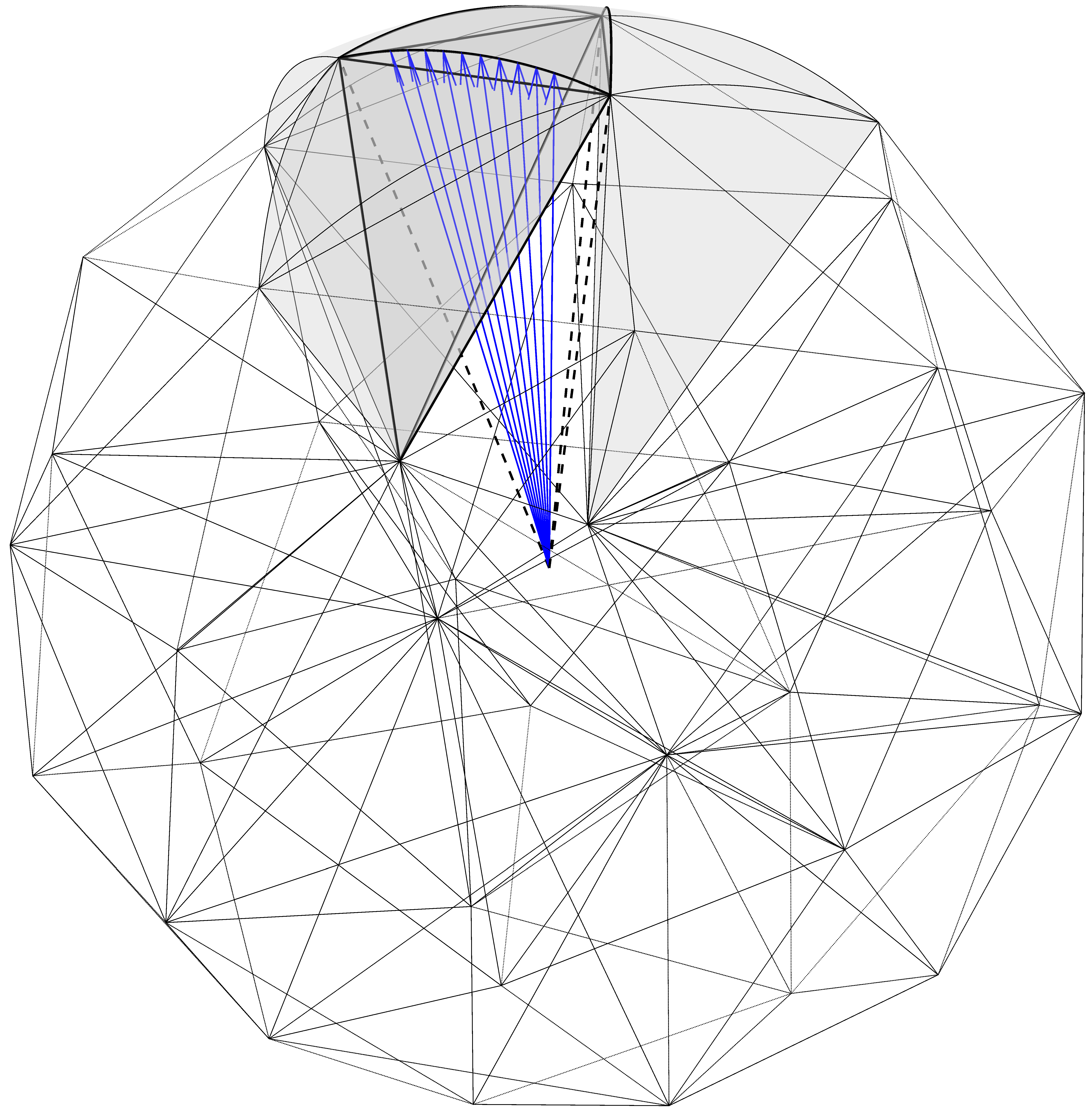}
\end{center}
\caption{An illustration of the tetrahedra in the set $T$.  The volume of one of the tetrahedra near the surface is outlined by thicker curves.  Call the outlined tetrahedron $t_{k}$ and let $\tau_{k,*}$ be the face of $t_{k}$ with three vertices on the surface of the sphere.  The arrows indicate the projection of one of $\tau_{k,*}$'s edges to the surface from the origin.  The dashed lines indicate the projection of $\tau_{k,*}$'s vertices from the origin.  When decomposing the volume of the sphere, the area of the triangle $\tau_{k,*}$, the areas between the arcs on the sphere and the edges of $\tau_{k,*}$ and the area of the spherical triangle between these arcs make of the boundary of the sliver of volume associated with $t_{k}$.  The projection ensures that between adjacent slivers there are no gaps or overlaps, illustrated by the slivers of volume associated with three adjacent tetrahedra.}
\label{fig:Sliver_Side_Projection}
\end{figure}

When assigning a sliver of volume to a particular tetrahedron, $t_{k}$, care must be taken so that there are no gaps or overlaps between adjacent slivers.  Let $\tau_{k,i}$, $i=1,2,3,4$, be the triangular faces of $t_{k}$.  At least one of these faces has all three vertices on the surface of the sphere.  In most cases, this will be only one face of $t_{k}$ (particularly when the volume is well resolved by small enough tetrahedra), call it $\tau_{k,*}$.  It turns out that if the three edges of $\tau_{k,*}$ are projected radially from the center of the sphere to the surface, gaps and overlaps will be prevented.  For each edge of $\tau_{k,*}$ the area between the arc on the sphere surface and the edge of the triangle forms a side of the sliver of volume.  The boundary of the sliver volume is formed by all three of these sides, the spherical triangle on the surface of the sphere between the three sides, and the triangle $\tau_{k,*}$.  Figure \ref{fig:Sliver_Side_Projection} illustrates one of these volumes.

Assigning the slivers of volume in this way provides for a transformation of the coordinates of the sliver which allows the integral over the volume to be written as an iterated integral over a triangular area and a parameter, $\sigma$, which relates to the projection from the origin.  Consider any point $\mathbf{x}$ in the volume.  The vector $\mathbf{x}$ intersects the plane containing $\tau_{k,*}$ at a point $\mathbf{x}'$ that is inside the triangle $\tau_{k,*}$.  All of the points inside $\tau_{k,*}$ can be parameterized by, for instance,
\begin{align}
    \mathbf{x}'(\lambda,\mu)=(1-\lambda)\mathbf{a}_{k}+\lambda\left((1-\mu)\mathbf{b}_{k}+\mu\mathbf{c}_{k}\right),\mbox{ }0\leq\lambda\leq 1\mbox{ and }0\leq \mu\leq 1,\nonumber
\end{align}
where $\mathbf{a}_{k}$, $\mathbf{b}_{k}$ and $\mathbf{c}_{k}$ are now representing the vertices of $t_{k,*}$.  With this parameterization of $t_{k,*}$ any point $\mathbf{x}$ in the volume of the sliver can be represented as
\begin{align}
    \mathbf{x}(\lambda,\mu,\sigma)=\left(1+\frac{\sigma}{\left\lVert\mathbf{x}'(\lambda,\mu)\right\rVert_{2}}\right)\mathbf{x}'(\lambda,\mu),\nonumber
\end{align}
where $0\leq \sigma \leq \rho-\left\lVert \mathbf{x}'(\lambda,\mu)\right\rVert_{2}$.   Here $\sigma$ measures the distance from $\mathbf{x}'$ to $\mathbf{x}$.  With this parameterization, for instance,
\begin{align}
    \iiint\limits_{s_{k}}\phi\left(\left\lVert \mathbf{x}-\mathbf{x}_{k,j}\right\rVert\right)dV=\int\limits_{0}^{1}\int\limits_{0}^{1}\int\limits_{0}^{\rho-\left\lVert \mathbf{x}'(\lambda,\mu)\right\rVert_{2}}\phi\left(\left\lVert \mathbf{x}(\lambda,\mu,\sigma)-\mathbf{x}_{k,j}\right\rVert\right) \left\lvert J(\sigma,\lambda,\mu)\right\rvert d\sigma d\lambda d\mu,\nonumber
\end{align}
where $J$ is the Jacobian determinant of $\mathbf{x}$ with respect to $\sigma$, $\lambda$ and $\mu$.  The integrals in $\sigma$, $\lambda$, and $\mu$ can be easily treated with any number of quadrature methods over an interval.

\subsection{Step 4: Combine weights from the subdomains}
Summing over all $k\in\{1,2,\ldots,K\}$ leads to the approximation of the volume integral over $\Omega$
\begin{align}
\iiint\limits_{\Omega}f(\mathbf{x})dV\approx\sum\limits_{k=1}^{K}\sum_{j=1}^{n}w_{k,j}f(\mathbf{x}_{k,j}). \nonumber
\end{align}

Let $\mathcal{K}_{i}$, $i=1,2,\ldots,N$, be the set of all pairs $(k,j)$ such that $\mathbf{x}_{k,j}\mapsto\mathbf{x}_{i}$. Then the volume integral over $\Omega$ can be rewritten as
\begin{align}
\iiint\limits_{\Omega}f(\mathbf{x})dV\approx\sum\limits_{i=1}^{N}W_{i}f(\mathbf{x}_{i}).\label{eq:quadrule}
\end{align}

\section{Test Examples} \label{sec:Test_Examples}

To demonstrate the performance of the method described herein, the algorithm will be applied to four different test integrands featuring varying degrees of smoothness.   The fourth of these test integrands includes extremely localized feature.   Weights are computed on quasi-uniformly spaced nodes, pseudo-randomly generated nodes, and a clustered node set with increased density near the localized feature of the fourth test integrand.  This clustered node set is used demonstrate the performance of the algorithm under localized node refinement.  In all of these tests, the radius of the ball is fixed to $\rho=\frac{6^{\frac{1}{3}}}{2\pi^{\frac{1}{3}}}$ so that the volume of the ball is equal to one.

\subsection{Node Sets}

\begin{figure}[h]
\begin{center}
\includegraphics[width=5in]{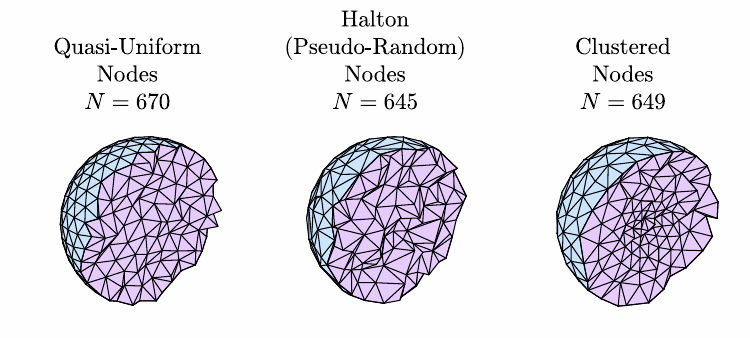}
\end{center}
\caption{Examples of the quasi-unifomly spaced, pseudo-randomly spaced, and clustered node sets.}
\label{fig:Ball_Node_Sets}
\end{figure}

In the cases of quasi-uniformly spaced nodes and the clustered node set, quadrature nodes were generated using a modification of the algorithm presented in \cite{PPGS04}.  The pseudo-randomly spaced node sets were generated by first drawing a set of points from the two-dimensional Halton sequence and mapping the set to the surface of the ball.  Then points were drawn from the three-dimensional Halton sequence, mapping the set to the domain $(x,y,z)\in[-\rho,\rho]\times[-\rho,\rho]\times[-\rho,\rho]$ and keeping only points satisfying $x^2+y^2+z^2\leq \rho-\frac{h}{10}$, where $h$ is prescribed to be the average spacing between the nodes on the surface.  Examples of the node sets are displayed in figure \ref{fig:Ball_Node_Sets}.

\subsection{Performance on Test Integrands}

\begin{figure}[h]
\begin{center}
\includegraphics[width=5in]{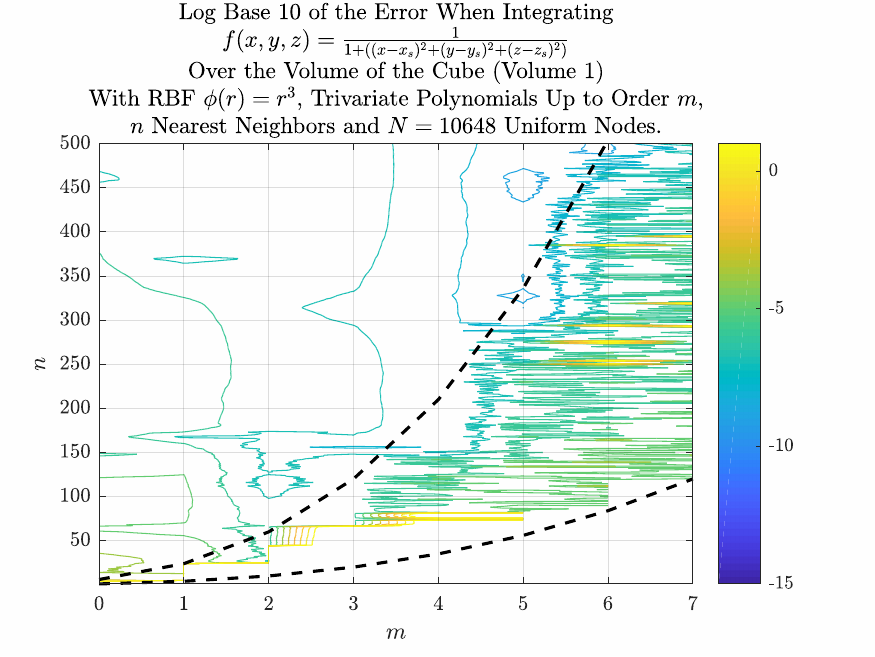}
\end{center}
\caption{Log base 10 of the absolute error when integrating $f(x,y,z)=\frac{1}{1+((x-x_{s})^{2}+(y-y_{s})^{2}+(z-z_{s})^{2})}$, with $x_{s}=0.234841098236337$, $y_{s}=0.048716273957102$ and $z_{s}=0.214415743035283$, over the volume of the unit cube centered at the origin for various choices of $n$ and $m$. The lower dashed curve is the function $n=\frac{(m+1)(m+2)(m+3)}{6}$ below which the matrix $A_{k}$ is guaranteed to be singular.  The upper dashed curve is $n=(m+1)(m+2)(m+3)$.  Choices of $n$ roughly above this curve lead to more accurate approximations of the integrand in the presence of boundaries.}
\label{fig:Cube_Errors_n_m}
\end{figure}

The algorithm was applied to four test integrands.  When generating quadrature weights, in all cases the radial basis function was $\phi(r)=r^3$ and the number of nearest neighbors, $n=(m+1)(m+2)(m+3)$, was based on the trivariate polynomial order $m$. Some computational experiments in, for instance, \cite{JARBF2017,VBNFBFGAB2017,VBNFBF2019} indicated that in the presence of boundaries the number of nearest neighbors must be large enough to overcome effects like Runge phenomenon.  The examples given in \cite{JARBF2017} indicated that the boundary errors were most prominent when nodes were (exactly) uniformly spaced.  Therefore, to determine how many nearest neighbors should be included to overcome boundary errors, the algorithm here was modified to compute volume integrals over a cube.  When considering a cube the entire volume can be decomposed by tetrahedra, so the algorithm need not consider slivers of volume near the surface. Figure \ref{fig:Cube_Errors_n_m} illustrates the absolute error when integrating $f(x,y,z)=\frac{1}{1+((x-x_{s})^{2}+(y-y_{s})^{2}+(z-z_{s})^{2})}$, with $x_{s}=0.234841098236337$, $y_{s}=0.048716273957102$ and $z_{s}=0.214415743035283$, over the volume of the unit cube centered at the origin for various choices of $n$ and $m$.  The matrix $A_{k}$ is singular for choices of $n$ below $n=\frac{(m+1)(m+2)(m+3)}{6}$, this is indicated by the lower dashed curve in the figure.  Further, for each case of $m=0,1,2,\ldots,7$ it is clear that choices of $n$ below $n=(m+1)(m+2)(m+3)$ can lead to large errors.

\begin{figure}[h]
\begin{center}
\includegraphics[width=5in]{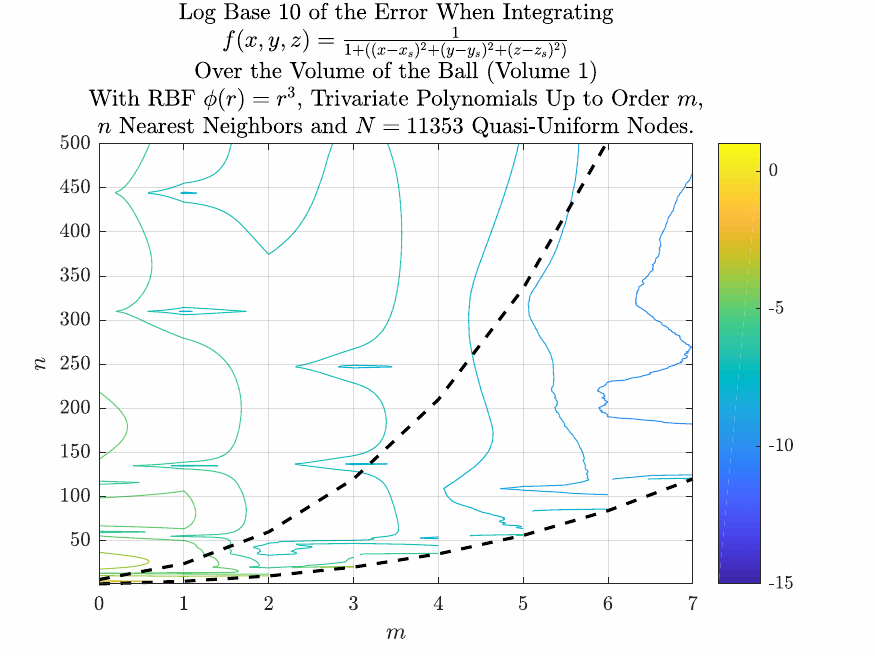}
\end{center}
\caption{Log base 10 of the absolute error when integrating $f(x,y,z)=\frac{1}{1+((x-x_{s})^{2}+(y-y_{s})^{2}+(z-z_{s})^{2})}$, with $x_{s}=0.234841098236337$, $y_{s}=0.048716273957102$ and $z_{s}=0.214415743035283$, over the volume of the unit ball centered at the origin for various choices of $n$ and $m$. The lower dashed curve is the function $n=\frac{(m+1)(m+2)(m+3)}{6}$ below which the matrix $A_{k}$ is guaranteed to be singular.  The upper dashed curve is $n=(m+1)(m+2)(m+3)$.}
\label{fig:Ball_Errors_n_m}
\end{figure}

Performing a similar experiment for the unit ball, figure \ref{fig:Ball_Errors_n_m} illustrates that even the relaxation from nodes that are exactly uniformly spaced to those that are quasi-uniformly spaces can allow for $n$ as low as $(m+1)(m+2)(m+3)$.  However, all results shown here utilize $n=(m+1)(m+2)(m+3)$.

\begin{figure}[h]
\begin{center}
\includegraphics[width=5in]{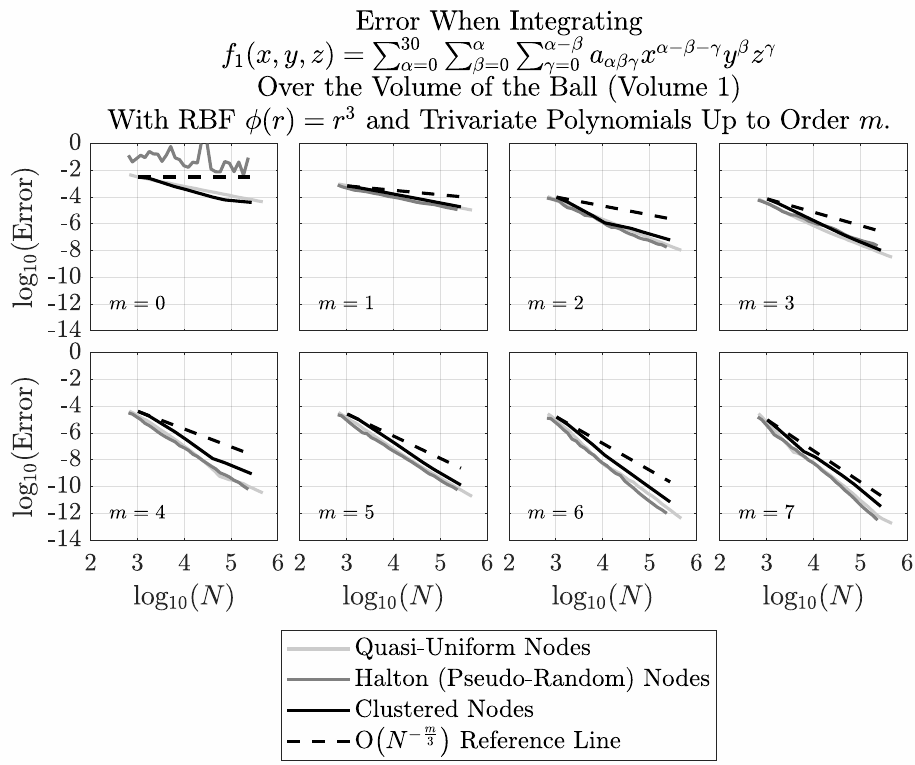}
\end{center}
\caption{Log base 10 of the absolute error when approximating the volume integral of $f_1$ over the ball with radius $\rho=\frac{6^{\frac{1}{3}}}{2\pi^{\frac{1}{3}}}$ centered at the origin.  The errors shown here are the largest after rotating the integrand 100 times.}
\label{fig:Ball_Quadrature_Error_f_1}
\end{figure}

The first of the test integrands is a degree 30 trivariate polynomial.  That is, let
\begin{align}
f_{1}(x,y,z)=\sum_{\alpha=0}^{30}\sum_{\beta=0}^{\alpha}\sum_{\gamma=0}^{\alpha-\beta}a_{\alpha\beta\gamma}x^{\alpha-\beta-\gamma}y^{\beta}z^{\gamma}.\nonumber
\end{align}
The coefficients of the polynomial are available in a Matlab file at \cite{BallVolumeQuadCodeMatlab}.  The exact value of the integral over the ball is
\begin{align}
\sum_{\alpha=0}^{30}\sum_{\beta=0}^{\alpha}\sum_{\gamma=0}^{\alpha-\beta}a_{\alpha\beta\gamma}\frac{\rho ^{\alpha +3} }{8 \Gamma \left(\frac{\alpha +5}{2}\right)}&\left[\left((-1)^{\beta }+1\right) \left((-1)^{\gamma }+1\right)\left(1+(-1)^{\alpha-\beta -\gamma }\right)  \right. \nonumber\\
&\left.  \Gamma \left(\frac{\beta +1}{2}\right) \Gamma \left(\frac{\gamma +1}{2}\right) \Gamma \left(\frac{\alpha -\beta -\gamma +1}{2}\right)\right]\nonumber
\end{align}
with $\Gamma$ the gamma function, and for the set of coefficients used here the expression evaluates to 3.792079311949332.  Figure \ref{fig:Ball_Quadrature_Error_f_1} displays convergence of the approximate integral to the exact value at an order better than O$\left(N^{-\frac{m}{3}}\right)$, where $m$ corresponds to the order of the trivariate polynomial terms used in the approximation.  If $h$ refers to a typical node separation distance, this corresponds to a convergence order of better than O$\left(h^{m}\right)$, especially in the case of quasi-uniformly spaced nodes. The theory in \cite{VB2019} explains that if the multivariate polynomial basis up to degree $m$ is included in the process of RBF interpolation, then all of the terms in the Taylor series up to degree $m$ will be handled exactly for the function being interpolated.  The remaining terms in the Taylor series are then approximated by the RBF basis that was included.  This is leading to a convergence order of at least O$(h^{m})$.

\begin{figure}[h]
\begin{center}
\includegraphics[width=5in]{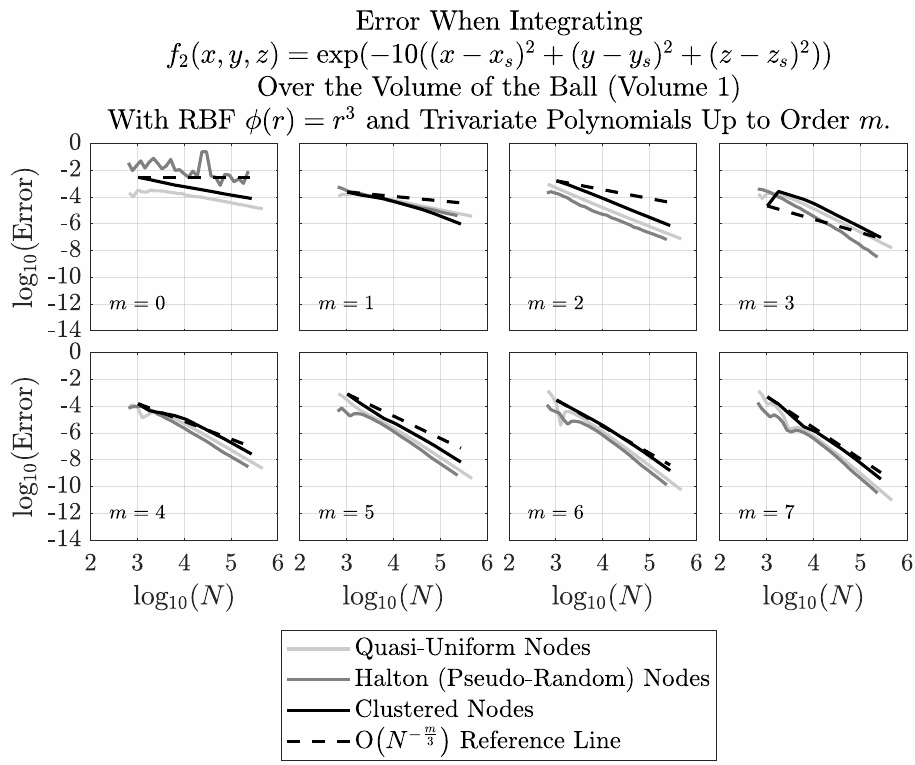}
\end{center}
\caption{Log base 10 of the absolute error when approximating the volume integral of $f_2$ over the ball with radius $\rho=\frac{6^{\frac{1}{3}}}{2\pi^{\frac{1}{3}}}$ centered at the origin.  The errors shown here are the largest after rotating the integrand 100 times.}
\label{fig:Ball_Quadrature_Error_f_2}
\end{figure}

The second test integrand is the Gaussian
\begin{align}
    f_{2}(x,y,z)=\exp\left(-10\left((x-x_{s})^{2}+(y-y_{s})^{2}+(z-z_{s})^2\right)\right)\nonumber
\end{align}
where
\begin{align}
(x_{s},y_{s},z_{s})=(0.047056440432708,0.071766893999009,0.118950756342700)\nonumber
\end{align}
is a randomly chosen shift of the center of the Gaussian from the origin. In order to have an accurate value to compare to, the volume integral was first approximated by evaluating
\begin{align}
\int\limits_{-\rho }^{\rho }\int\limits_{-\sqrt{\rho ^2-x^2}}^{\sqrt{\rho ^2-x^2}}\int\limits_{-\sqrt{\rho ^2-x^2-y^2}}^{\sqrt{\rho ^2-x^2-y^2}}f_{2}(x,y,z)dzdydx,\nonumber
\end{align}
using Matlab's integral3 command with the absolute and relative tolerances both set to ten times machine precision.  The resulting approximation of the integral for $f_{2}$ is 0.161965667295343.  Figure \ref{fig:Ball_Quadrature_Error_f_2} illustrates the error in the integral of $f_{2}$ over the ball when compared to the result from Matlab after rotating the integrand randomly 100 times.  It is clear again that the order of the error is most dependent on the degree of the polynomials used in the interpolation.

\begin{figure}[h]
\begin{center}
\includegraphics[width=5in]{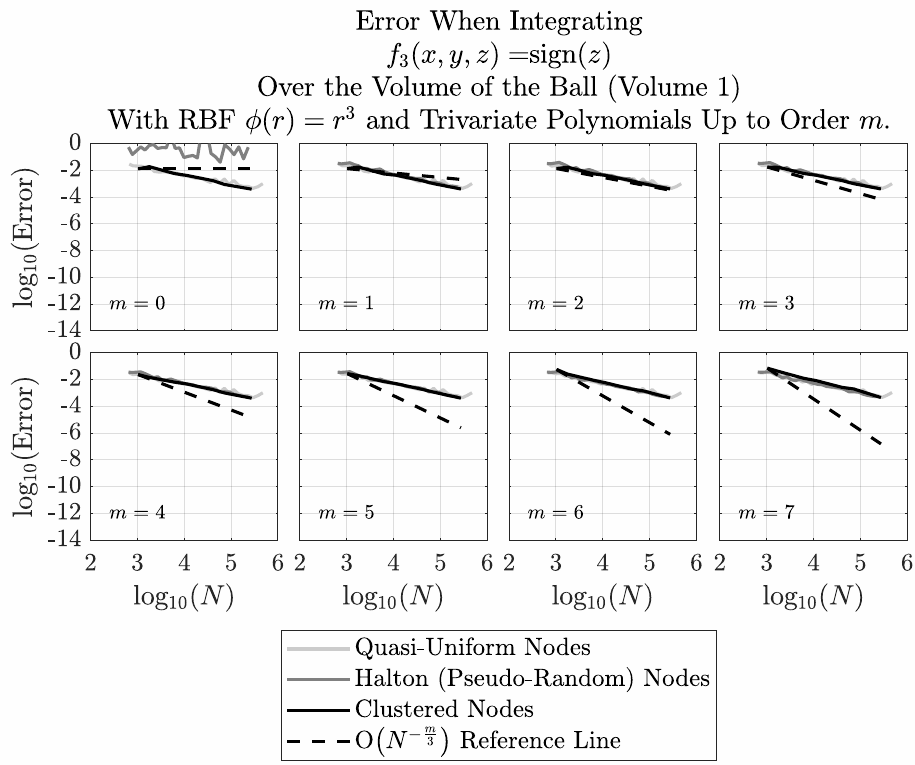}
\end{center}
\caption{Log base 10 of the absolute error when approximating the volume integral of $f_3$ over the ball with radius $\rho=\frac{6^{\frac{1}{3}}}{2\pi^{\frac{1}{3}}}$ centered at the origin.  The errors shown here are the largest after rotating the integrand 100 times.}
\label{fig:Ball_Quadrature_Error_f_3}
\end{figure}

The third test integrand is
\begin{align}
    f_{3}(x,y,z)=\mbox{sign}\left(z\right)\nonumber
\end{align}
with $\mbox{sign}$ the signum function.  This function is discontinuous at the plane $z=0$, so any method based on a continuous approximation of the integrand across the discontinuity should not be expected to achieve better than O$\left(N^{-\frac{1}{3}}\right)$ (i.e. O$(h)$) error.  Figure \ref{fig:Ball_Quadrature_Error_f_3} illustrates this for the present method.

\subsection{Performance When Utilizing Clustered Node Sets}

\begin{figure}[h]
\begin{center}
\includegraphics[width=5in]{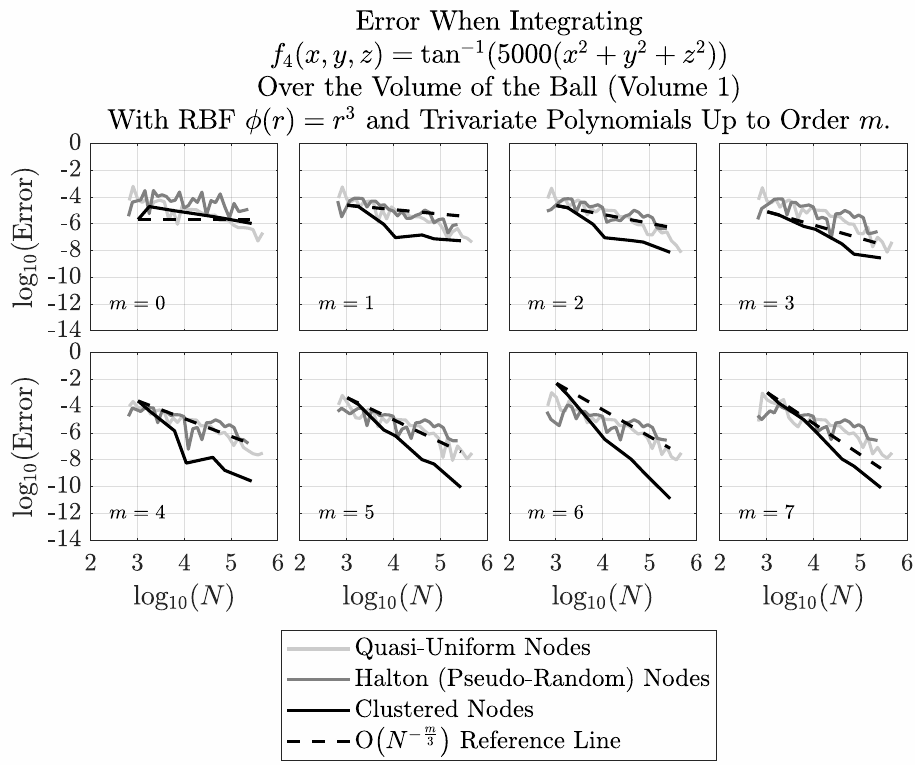}
\end{center}
\caption{Log base 10 of the absolute error when approximating the volume integral of $f_4$ over the ball with radius $\rho=\frac{6^{\frac{1}{3}}}{2\pi^{\frac{1}{3}}}$ centered at the origin. }
\label{fig:Ball_Quadrature_Error_f_4}
\end{figure}

To illustrate further utility of the proposed method, the algorithm was also applied to a test integrand featuring a steep and localized gradient.  To capture the rapid change in the integrand node sets were generated that feature more densely clustered nodes near the local feature.  The test integrand was
\begin{align}
    f_{4}(x,y,z)=\tan^{-1}\left(5000(x^2+y^2+z^2)\right),\nonumber
\end{align}
which has a steep gradient near the origin.  Figure \ref{fig:Ball_Quadrature_Error_f_4} illustrates that in cases where quasi-uniformly spaced or pseudo-randomly spaced node sets cannot capture the changes in the integrand, weights generated for node sets with clustering near features improve the approximation under refinement.

\subsection{Computational Expense}

\begin{figure}[h]
\begin{center}
\includegraphics[width=4.1in]{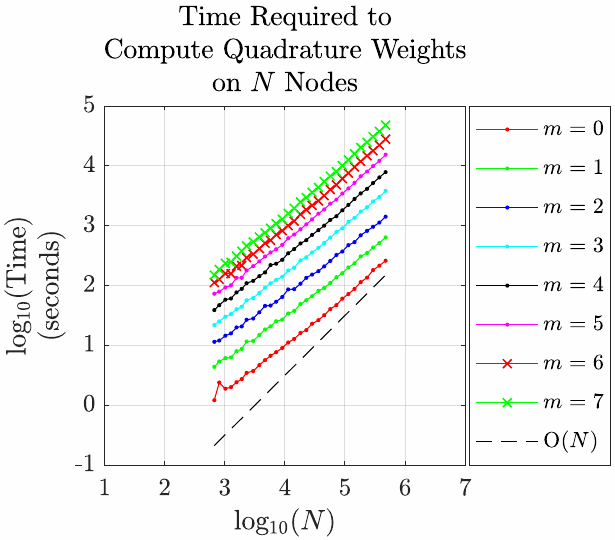}
\end{center}
\caption{Log base 10 of the time it takes to compute quadrature weights on $N$ nodes when including trivariate polynomial up to order $k$.  The (black) dashed line is an O$(N)$ reference line.}
\label{fig:Ball_Times}
\end{figure}

Just like the algorithms presented in \cite{JARBF2016,JARBFMLW2016,JARBF2017}, the ability to consider each tetrahedron individually allows the time to compute a set of quadrature weights and the use of memory both to scale like $O(N)$. Figure \ref{fig:Ball_Times} illustrates the time to compute the set of quadrature weights on $N$ nodes for various choice of the polynomial order, $m$.  Since the choice of $m$ affects the sizes of the systems of linear equations that need to be solved at each iteration, the figure shows an increase in the computational cost as $m$ increases.
Further, except for the identification of nearest neighbors in order to construct the local weight set for each tetrahedron/sliver of volume and for the combination of weights in step 4 the algorithm is pleasingly parallel.  The parallelization tests in \cite{JARBF2016} illustrate that these two steps do not have a significant impact on the scalability of the algorithm with the number of cores when considering evaluating surface integrals, and the same is true for the volume integration algorithm described herein.

\section{Conclusions} \label{sec:Conclusions}

This study has supplemented the previous RBF-FD based approach for evaluating definite integrals \cite{JARBF2016,JARBFMLW2016,JARBF2017} with an extension to integrals over volumes. The computational tests illustrate an algorithm that can achieve at least $O(h^{m})$ accuracy, with $h$ the typical node separation distance and $m$ the order of trivariate polynomial basis functions included in the approximation.  On a set of $N$ nodes in the ball, the computational cost is only $O(N)$ and the algorithm is pleasingly parallel.  A key feature of the algorithm is that it is able to compute quadrature weights on even irregularly spaced or clustered node sets.


\bibliographystyle{siamplain}
\bibliography{references}
\end{document}